\documentclass[10pt]{article}
\usepackage{times,amsfonts,amsmath,amssymb,exscale,mathrsfs}
\usepackage{pst-all}
\usepackage{graphicx}
\usepackage{xcolor} 
\usepackage{epsfig}
\usepackage{shuffle}
\usepackage[numbers]{natbib}
\newcommand{\C}{\mathbb{C}}
\newcommand{\R}{\mathbb{R}}
\newcommand{\wht}{\widehat t}
\newcommand{\W}{\mathscr{W}}
\newcommand{\Z}{\mathbb{Z}}

\begin{document}

\title{Word-series high-order averaging of highly oscillatory  differential equations with delay
}
\author
{J. M. Sanz-Serna \\
\footnotesize Departamento de Matem\'{a}ticas, Universidad Carlos III de Madrid\\
\footnotesize Avenida de la Universidad 30, E-28911 Legan\'{e}s (Madrid), Spain\\
\footnotesize Email: jmsanzserna@gmail.com
\\[2mm]
Beibei Zhu\\
\footnotesize National Center for Mathematics and Interdisciplinary Science\\
\footnotesize Academy of Mathematics and Systems Science\\
\footnotesize Chinese Academy of Sciences, Beijing 100190, P.R. China\\
\footnotesize Email: zhubeibei@lsec.cc.ac.cn
}

\date{}
\maketitle

\begin{abstract}
We show that, for appropriate combinations of the values of the delay and the forcing frequency, it is possible to obtain easily high-order averaged versions of periodically forced systems of delay differential equations with constant delay. Our approach is based on the use of word-series techniques to obtain high-order averaged equations  for differential equations without delay.
\end{abstract}
\bigskip

{\footnotesize\noindent\textbf{Mathematical Subject Classification (2010)} 34C29

}

\section{Introduction}
We show that, for appropriate combinations of the values of the delay and the forcing frequency, it is possible to obtain easily high-order averaged systems of periodically forced systems of delay differential equations with constant delay.

The well-known theory of averaging  \cite{sanders,lehman} studies the reduction, by means of time-dependent changes of variables, of (nonautonomous) forced systems of differential equations to  autonomous time-independent systems (averaged systems). Such a reduction is useful because autonomous systems are  easier to analyze than their nonautonomous counterparts. In addition, the numerical integration of periodically or quasi-periodically forced systems may be expensive because  typically  integrators have to operate with step-sizes that are small with respect to the shortest period present in the forcing; in those cases, integrating an averaged version of the given system may be very advantageous.

In many occasions averaging is only carried out to first order. When the forcing is periodic this simply means replacing the right-hand side of the given oscillatory differential system by its average over one period of the forcing. In other cases, first order averaging is not enough and it is desirable to obtain averaged systems that provide more accurate approximations. For ordinary differential equations without delay the obtention of high-order averaging may be a daunting task, even with the help of a symbolic manipulator. The series of papers \cite{part1,part2,orlando,part3,guirao,erratum,icmat, abel} has developed a technique that, in the absence of delays, makes it possible to compute high-order averaged systems through simple recursions. The technique is based on the  theory of \emph{word series} \cite{focm16}, formal series \cite{china} that have numerous applications in the fields of deterministic or stochastic dynamical systems \cite{juanluis} and numerical integration \cite{alfonso,alfonso2}.

The difficulties in obtaining high-order averaged systems are compounded if the system to be averaged has \emph{delays}. In this paper we show that, for periodically forced differential systems with constant delay, it is possible to obtain high-order averaged systems by an application of the word-series results in \cite{part1,part2,orlando,part3,guirao,erratum,icmat, abel}. The simple treatment presented here is only possible when the forcing period is a submultiple of the delay, a hypothesis whose scope is discussed in Section 3.

Section 2 recalls the  word-series averaging of systems of differential equations without delay. Section 3 contains the main idea and Section 4 describes an application that shows how third-order averaging succeeds where second-order averaging fails. The final Section 5 concludes.
\section{Preliminaries}
In this section we summarize the word-series averaging of differential equations without delay \cite{icmat}. This summary will provide the foundation to address later the case with delay.

Let us consider highly oscillatory initial value problems of the form
\begin{equation}\label{eq:oscode}
\frac{d}{dt}\xi = g(\xi, \Omega t),\quad t\geq0,\qquad \xi(0) = \xi_0\in\C^D.
\end{equation}
Here the smooth function \(g(\xi,\theta)\) is \(2\pi\)-periodic in its second argument \(\theta\in\R\), with Fourier expansion
\begin{equation}\label{eq:fourier}
g(\xi,\theta) = \sum_{k\in\Z}  \exp(ik\theta) g_k(\xi),
\end{equation}
and the parameter \(\Omega\gg 1\) is the angular frequency of the fast periodic forcing. The corresponding period is
\(T=2\pi/\Omega\).

In the theory of word series the set of indices \(\Z\) in \eqref{eq:fourier} is seen as an (infinite) alphabet and for \(n=1,2,\dots\), the elements of \(\Z^n\) are regarded as words consisting of \(n\) letters of this alphabet. There is also an empty word \(\emptyset\) with \(n=0\) letters. Given \eqref{eq:oscode}, to each word
\(w\)  one associates the corresponding \emph{word basis function} \(g_w\), a map \(\C^D\rightarrow \C^D\). By definition, for words with one letter \(k\in\Z\), the word basis function \(g_w(\xi)\) is the Fourier coefficient \(g_k(\xi)\) in \eqref{eq:fourier}. For words with \(n>1\) letters, we define recursively
\begin{equation}\label{eq:wordbasis}
g_{k_1\cdots k_n}(\xi) = g^\prime_{k_2\cdots k_n}(\xi)g_{k_1}(\xi),
\end{equation}
where \(g^\prime_{k_2\cdots k_n}(\xi)\) is the Jacobian matrix of \(g_{k_2\cdots k_n}(\xi)\). For the empty word \(\emptyset\), the basis function is the identity map \(\xi\mapsto\xi\).

The set of all words (including the empty word) is denoted by \(\W\) and \(\C^\W\) refers to the set of all mappings \(\delta:\W\rightarrow \C\). For \(\delta\in\C^\W\) and \(w\in\W\), \(\delta_w\) denotes the value of \(\delta\) at \(w\). To each \(\delta\in\C^\W\) we associate the corresponding \emph{word series} (relative to the mappings \(g_k\) in \eqref{eq:fourier}); this is the formal series
\[
W_\delta(\xi) = \sum_{w\in\W}\delta_w g_w(\xi),
\]
whose terms are maps \(\C^D\rightarrow \C^D\). The complex numbers \(\delta_w\) are  the \emph{coefficients} of the series. As an example consider the case where \(\delta_\emptyset = 1\) and \(\delta_w = 0\) for each nonempty word \(w\); in this case \(W_\delta(\xi)=\xi\).

As proved in \cite{part2,icmat} (see \cite{abel} for an alternative technique),  there exist  \(\bar\beta\in\C^\W\) and a \(2\pi\)-periodic map \(\theta\in\R\rightarrow \kappa(\theta)\in \C^\W\)   such that the solution  of \eqref{eq:oscode} may be written as
\begin{equation}\label{eq:change}
\xi(t) = W_{\kappa(\Omega t)}(\Xi(t)),
\end{equation}
where \(\Xi(t)\) is the solution of
\begin{equation}\label{eq:averode}
\frac{d}{dt}\Xi = W_{\bar \beta}(\Xi),\qquad \Xi(0) = \xi_0\in\R^D.
\end{equation}
In this way, the time-dependent change of variables \(\xi\mapsto \Xi\) in \eqref{eq:change} transforms the given highly oscillatory initial value problem into the initial value problem \eqref{eq:averode} where there is no periodic forcing. Therefore \eqref{eq:averode} provides an \emph{averaged} version of \eqref{eq:oscode}. In addition,
the coefficients \(\kappa_w(\theta)\) are such that, at \(\theta = 2k\pi\), \(k=0,1,2,\dots\), we have \(\kappa_\emptyset(\theta)=1\) and \(\kappa_w(\theta)= 0\) for each nonempty word \(w \). This implies that, for those values of \(\theta\), the transformation \(\xi \mapsto \Xi\) is the identity map and it then follows that,
in \eqref{eq:change}, \(\xi(t) = \Xi(t)\) at the \emph{strobocopic times}, \(t = 0,T, 2T, \dots\) We then say that the change of variables is stroboscopic and that \eqref{eq:averode} is obtained from \eqref{eq:oscode} through \emph{stroboscopic averaging}. There are of course alternative forms of averaging; for instance one may impose the condition that the periodic change of variables \(\xi\mapsto\Xi\) has 0 average over one period, rather than being the identity at \(\theta = 0\).

It is in order to point out that \(\bar \beta\) and \(\kappa(\theta)\) depend on \(\Omega\), but are otherwise \emph{universal} in the sense that they do not change with the dimension \(D\) or with the choice of \(g(\xi,\theta)\). Such a universality implies that they may be computed once and for all; averaging a new differential system requires the computation of new basis functions but not of new coeffcients. The values \(\bar\beta_w\) and \(\kappa_w(\theta)\) for \(w\in\W\) may be computed easily by recursion with respect to the number of letters in the word \(w\); the reader is referred to \cite{part2,erratum,icmat} for details.

In general, the formal series  \(W_{\bar\beta}(\Xi)\) in the averaged system \eqref{eq:averode} and the formal series in the change of variables \eqref{eq:change} do not converge and have to be truncated and of course the truncation introduces an error, see \cite{orlando,part3}. For  our purposes here we just mention that, if \(W_{\bar\beta}(\Xi)\) is truncated  by eliminating all the terms in the series that correspond to words with more than \(n\) letters, one obtains an initial value problem whose solution \(\Xi(t)\) coincides at stroboscopic times with the solution  \(\xi(t)\) except for an error of size \(\mathcal{O}(1/\Omega^n)\) as \(\Omega \rightarrow \infty\).\footnote{To approximate \(\xi\) at times that are not stroboscopic one needs to apply to the solution of the truncated averaged system a change of variables obtained by truncating the series in \eqref{eq:change} .} The truncated averaged system with first order errors obtained by discarding contributions corresponding to words with two or more letters  is found to be:
\begin{equation}\label{eq:first}
\frac{d}{dt} \Xi = g_0(\Xi),
\end{equation}
as it may have been expected.

The truncated averaged system with second order errors \(\mathcal{O}(1/\Omega^2)\) is given by \cite{part2,erratum,icmat}
\begin{equation}\label{eq:second}
\frac{d}{dt} \Xi = g_0(\Xi)+\sum_{k\neq 0} \frac{i}{k\Omega}
\Big(
g_0^\prime(\Xi)g_k(\Xi)-
g_k^\prime(\Xi)g_0(\Xi)
+g_k^\prime(\Xi)g_{-k}(\Xi)
\Big).
\end{equation}
Note that, \(g_0\) is by definition the word-basis function associated with the one-letter words \(0\). Furthermore, according to \eqref{eq:wordbasis}, the function \(g_0^\prime g_k\) is the word basis function associated to the word \(k0\) and similarly \(g_k^\prime g_0\), \(g_k^\prime g_{-k}\) are associated to \(0k\) and \(-kk\). Thus the right-hand side of \eqref{eq:second} is indeed a word series (whose coefficients vanish for all words with 3 or more letters).

 The closed form expression for the third order averaged system for the problem \eqref{eq:oscode} with general \(g\) may be seen in \cite{part2,erratum,icmat}. For order \(\geq 4\) it is not practical to  find the general expression of the averaged system and then to apply it to the specific instance of \eqref{eq:oscode} at hand. One should rather find recursively  the word basis functions corresponding to the \(g\) of interest, perhaps with the help of a symbolic manipulator, see \cite{guirao,abel} for additional details.

\section{Highly oscillatory problems with delay}
We now consider the \(D\)-dimensional system with constant delay \(\tau>0\)
\begin{eqnarray}\label{eq:osc}
\frac{d}{dt}x(t)&=&f(x(t),x(t-\tau),\Omega t),\qquad t\ge 0\\
x(t)&=&\varphi(t),\qquad -\tau\le t\le 0,\label{eq:osc2}
\end{eqnarray}
where \(f(x,y,\theta)\) is \(2\pi\)-periodic in its third argument, with Fourier expansion
\begin{equation}\label{eq:fourier2}
f(x,y,\theta) = \sum_{k\in\Z}  \exp(ik\theta)f_k(x,y),
\end{equation}
and \(\Omega\) is the angular frequency as above. Without loss of generality \cite{beibei}, it is assumed that the known function
\(\varphi\) that specifies the initial history is independent of \(\Omega\). Reference \cite{beibei} contains a list of references to problems of this form that appear in different applications.

The developments that follow are based on the \emph{standing hypothesis} that the delay \(\tau\) is an integer multiple of the period \(T\), or in other words that \(t=\tau\) \emph{is a stroboscopic time}. Let us discuss this hypothesis. In some applications, the value of \(\tau\) is given and the interest is in studying the behaviour of the solution of \eqref{eq:osc}--\eqref{eq:osc2} for large \(\Omega\), but there is some freedom as to the choice of the specific value of \(\Omega\). In those cases, the theory below applies by choosing \(\Omega\) in such a way that \(\tau\Omega/(2\pi)\) is an integer. Similarly, our theory applies to situations where the interest is in a given, large value of \(\Omega\) and there is some freedom in the choice of \(\tau\) to be used in the model. On the other hand, if the values of \(\tau\) and \(\Omega\) are fixed, the material in this paper does not apply. In those cases finding high-order averaged systems may be extremely complicated (see e.g. the computations in \cite{beibei}).

\subsection{The averaging procedure}
\begin{figure}[t]
\vspace{-4cm}\centering\includegraphics[scale=0.45]{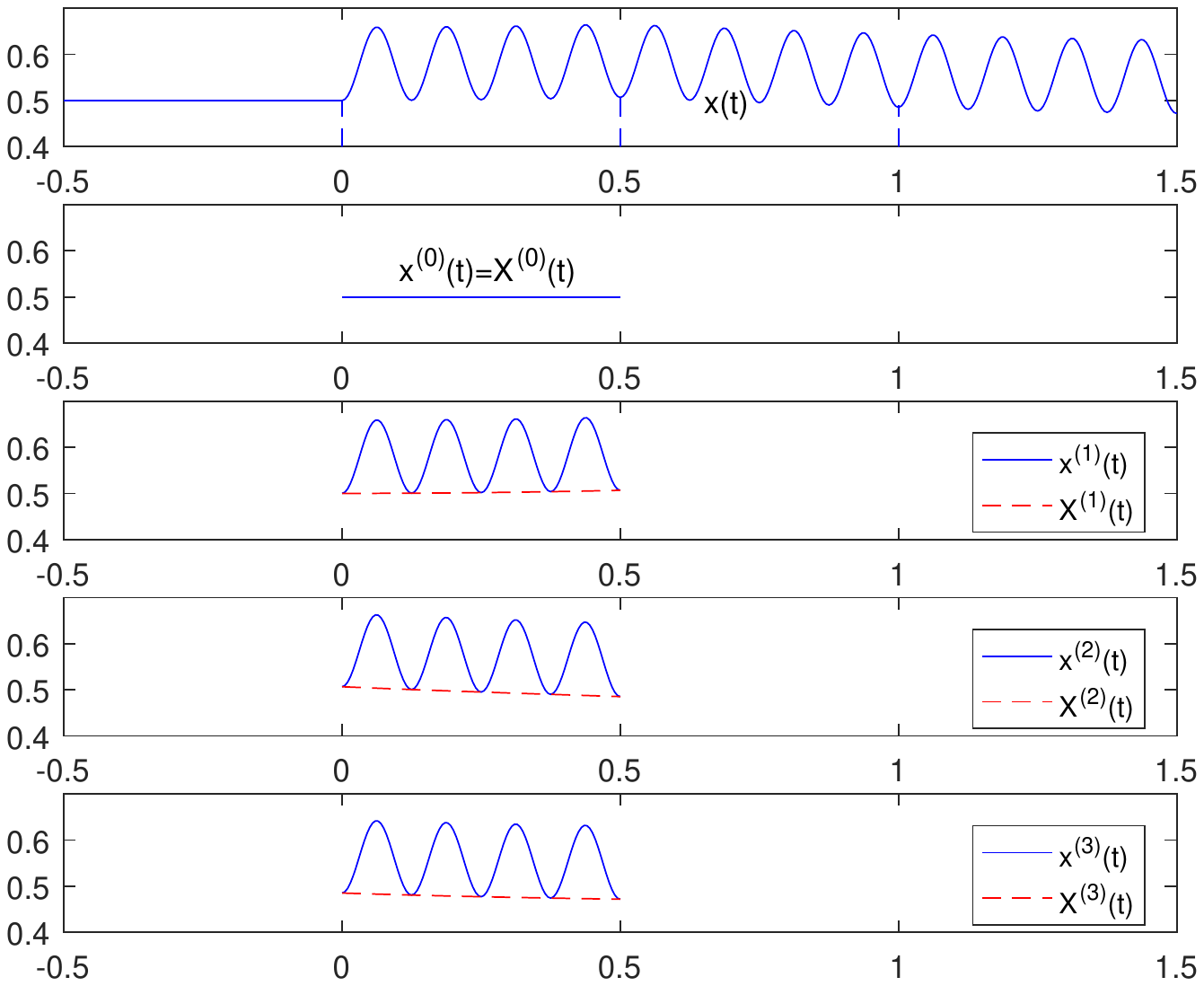}
\vspace{-4cm}
\caption{\footnotesize The solution \(x\) of the delay oscillatory problem (top subplot) is represented by an array \((x^{(0)},\dots,x^{(L)})\) of functions defined in the interval \(0\leq t\leq \tau\); the four bottom subplots depict these in the case \(L=3\). Averaging the ordinary differential system without delay satisfied by the \((x^{(0)},\dots,x^{(L)})\) leads to a differential system without delay for the averaged solution \((X^{(0)},\dots,X^{(L)})\). The functions \(X^{(\ell)}\) (dotted lines) are patched together to get the solution  \(X\) of the averaged delay problem that we wish to find.
}
\label{fig:B}
\end{figure}
Assume that \eqref{eq:osc}--\eqref{eq:osc2} is to be studied in a bounded interval of the form \(0\leq t\leq L\tau\) for a suitable integer \(L>0\). (This hypothesis is made at this stage for mathematical convenience to avoid systems of infinitely many differential equations; the averaged systems that we will find will be valid for all \(t\geq 0\).) We introduce the functions (see Fig.~\ref{fig:B}).
\begin{eqnarray}
x^{(0)}(t)&=&\varphi(t-\tau), \qquad 0\leq t \leq \tau, \label{eq:def1}\\
x^{(\ell)}(t)&=&x(t+(\ell-1)\tau), \qquad 0\leq t \leq \tau,\qquad \ell=1, \ldots, L,\label{eq:def2}
\end{eqnarray}
and note that determining these functions is clearly equivalent to determining the solution of \eqref{eq:osc}--\eqref{eq:osc2} in the interval \(0\leq t\leq L\tau\). The \(x^{(\ell)}(t)\) with \(\ell>0\) satisfy
the conditions
\begin{equation}
\label{eq:ode2}
x^{(\ell)}(0)=x^{(\ell-1)}(\tau),\qquad  1 \leq \ell \leq L.
\end{equation}
and the differential equations
\[
\frac{d}{dt}x^{(\ell)}(t)=f(x^{(\ell)}(t),x^{(\ell-1)}(t),\Omega (t+(\ell-1)\tau)),\quad
0\leq t \leq \tau,
\: 1 \leq \ell \leq L.
\]
Even though the relations \eqref{eq:ode2} may be reminiscent of a two-point boundary value problem, we are facing here an initial value problem: \(x^{(1)}(t)\) is determined by solving the differential equation with \(\ell =1\) with initial condition \(x^{(1)}(0)= \varphi(0)\), once \(x^{(1)}(t)\) is known, \(x^{(2)}(t)\) is determined by solving the differential equation with \(\ell =2\) with initial condition \(x^{(2)}(0)= x^{(1)}(\tau)\), etc.

By taking into account the periodicity of \(f(x,y,\theta)\) as a function of  \(\theta\) and our standing hypothesis that \(\tau\) is a multiple of the period, the last display may be written as
\begin{equation}\label{eq:system}
\frac{d}{dt}x^{(\ell)}(t)= f(x^{(\ell)}(t),x^{(\ell-1)}(t),\Omega t),\quad
0\leq t \leq \tau,
\: 1 \leq \ell \leq L.
\end{equation}
Note that this would be a system of the form \eqref{eq:oscode} for the \(D\times L\)-dimensional vector \((x^{(1)},\dots, x^{(L)})\) if it were not for the fact that the right-hand side of the equation corresponding to \(\ell=1\) contains \(x^{(0)}\), a known function of \(t\) given in \eqref{eq:def1}. In order to have a system of the form \eqref{eq:oscode}, where the right-hand side only depends on \(t\) through the combination \(\theta=\Omega t\), we
proceed as follows. We introduce the \(1+D(L+1)\)-dimensional vector of unknown functions of the variable \(t\)
\[
\xi=(\wht, x^{(0)},x^{(1)}, \dots, x^{(L)}),
\]
and add to \eqref{eq:system} the differential equations
\begin{equation}\label{eq:system2}
\frac{d}{dt} \wht = 1,\qquad \frac{d}{dt} x^{(0)} = \dot \varphi(\wht-\tau),\quad 0\leq t\leq \tau,
\end{equation}
(the dot represents differentiation) subject to the initial conditions
\begin{equation}\label{eq:ic}
\wht(0) = 0,\qquad x^{(0)}(0) = \varphi(-\tau).
\end{equation}
The solution of the initial value problem \eqref{eq:system2}--\eqref{eq:ic} is obviously \(\wht = t\) and \(x^{(0)}(t) = \varphi(t-\tau)\).
Now \eqref{eq:system2} in tandem with \eqref{eq:system} is a \(1+D(L+1)\)-dimensional system of the form \eqref{eq:oscode}
that may be stroboscopically averaged by following the procedure described in the previous section. If \(X^{(\ell)}\) is the averaged counterpart of
\(x^{(\ell)}\), we are interested in the solutions of the averaged system that satisfy (cf.~\eqref{eq:ode2})
\begin{equation}\label{eq:conditionsX}
X^{(\ell)}(0)=X^{(\ell-1)}(\tau),\qquad  1 \leq \ell \leq L,
\end{equation}
so that it is possible to define a continuous function \(X\) in the interval \([-\tau,L\tau]\) by patching together the different \(X^{(\ell)}\),
\begin{equation}\label{eq:patching}
X(t) = X^{(\ell)}(t-(\ell-1)\tau),\quad  -\tau\leq t \leq L\tau, \quad \ell = 0, \dots, L,
\end{equation}
thus undoing the process that we used to move from the solution \(x\) of \eqref{eq:osc}--\eqref{eq:osc2} to the functions \(x^{(\ell)}\) in \eqref{eq:def1}--\eqref{eq:def2}, see Fig.~\ref{fig:B}. In this way the averaged delay equation is obtained by patching the equations for the \(X^{(\ell)}\). The whole procedure will be clear after we present the simple case of the first-order averaged system.

\subsection{The first-order averaged system}
If \(g\) denotes the right-hand side of the oscillatory system of differential equations \eqref{eq:system}--\eqref{eq:system2} satisfied by \(\xi\), the Fourier series for \(g(\xi,\theta)\) (cf.~\eqref{eq:fourier}) has coefficients
\[
g_0 =
\left[
\begin{matrix}
1\\ \dot \varphi(\wht-\tau)\\ f_0(x^{(1)},x^{(0)})\\f_0(x^{(2)},x^{(1)})\\ \vdots\\f_0(x^{(L)},x^{(L-1)})
\end{matrix}
\right],\qquad
g_k =
\left[
\begin{matrix}
0\\ 0\\ f_k(x^{(1)},x^{(0)})\\f_k(x^{(2)},x^{(1)})\\ \vdots\\f_k(x^{(L)},x^{(L-1)})
\end{matrix}
\right],\quad k\neq 0,
\]
where the \(f_k\) are the Fourier coefficients of \(f\) as in \eqref{eq:fourier2}. According to \eqref{eq:first}, the first-order averaged system of differential equation without delay in the interval \(0\leq t\leq \tau\) is therefore given by (capital letters denote averaged dependent variables)
\begin{eqnarray*}
\frac{d}{dt} \widehat T &=& 1,\\
\frac{d}{dt} X^{(0)} & = & \dot \varphi(\widehat T-\tau),\\
\frac{d}{dt} X^{(\ell)} & = & f_0(X^{(\ell)}, X^{(\ell-1)}),\quad \ell = 1,\dots, L.
\end{eqnarray*}
This system is to be solved with the conditions \(\widehat T(0) = 0\), \(X^{(0)}(0) = \varphi(-\tau)\) and \eqref{eq:conditionsX}. Clearly
\(\widehat T = t\), \(X^{(0)}(t)=\varphi(t-\tau)\) and therefore the function \(X\) in \eqref{eq:patching} satisfies the averaged delay problem
(cf.~\eqref{eq:osc}--\eqref{eq:osc2})
\begin{eqnarray}\label{eq:aver1}
\frac{d}{dt}X(t)&=&f_0(X(t),X(t-\tau)),\\
X(t)&=&\varphi(t),\qquad -\tau\le t\le 0,\label{eq:aver12}
\end{eqnarray}
as one could have easily guessed. Note that \(L\) does not feature in the averaged problem, as we pointed out above.

The solutions \(x\) and \(X\) of \eqref{eq:osc}--\eqref{eq:osc2} and \eqref{eq:aver1}--\eqref{eq:aver12} differ at stroboscopic times  by
\(\mathcal{O}(1/\Omega)\) as we prove next. That \(x(t)-X(t) = \mathcal{O}(1/\Omega)\) at stroboscopic times \(t\leq \tau\) follows by considering that, for those values of \(t\),
\(x(t) = x^{(1)}(t)\), \(X(t) = X^{(1)}(t)\) and  \(x^{(1)}(t)-X^{(1)}(t) = \mathcal{O}(1/\Omega)\) because  of the properties of stroboscopic averaging of ordinary differential equations without delay. For stroboscopic times \(t\) in the interval \(\tau<t\leq 2\tau\) the situation is slightly more complicated because in addition to the error \(\mathcal{O}(1/\Omega)\) arising from truncating \eqref{eq:averode} there is an additional source of error: in this interval
 \(x(t) = x^{(2)}(t+\tau)\), \(X(t) = X^{(2)}(t+\tau)\) but \(x^{(2)}(t)\) and
\(X^{(2)}(t)\) do not share the same initial value at \(t=0\). However the difference of initial values is \(x^{(2)}(0)-X^{(2)}(0) =
x^{(1)}(\tau)-X^{(1)}(\tau)\), which we know is of size \(\mathcal{O}(1/\Omega)\). It follows that \(x(t)-X(t) = \mathcal{O}(1/\Omega)\) also at
stroboscopic times \(\tau<t\leq 2\tau\). The iteration of this argument leads to the conclusion we seek, i.e. \(x(t)-X(t) = \mathcal{O}(1/\Omega)\)
whenever \(t\) is a stroboscopic time. Of course the constant implied in the \(\mathcal{O}\) notation of course depends on \(t\) (and in general will grow as \(t\) increases).

\subsection{Second-order averaged system}
If rather than using the truncated averaged system \eqref{eq:first}, we use \eqref{eq:second}, the procedure described above leads, after some simple algebra, to
an averaged problem given by \eqref{eq:aver12} and a delay differential equation of the form
\begin{eqnarray}\label{eq:f21}
\frac{d}{dt}X(t) &=& F_{2,1}(X(t),X(t-\tau)), \quad 0\leq t <\tau,\\ \
\frac{d}{dt}X(t) &=& F_{2,2}(X(t),X(t-\tau),X(t-2\tau)), \quad t\geq \tau,
\label{eq:f22}
\end{eqnarray}
with
\begin{eqnarray*}
F_{2,1}&=&f_0(X(t),X(t-\tau))\\
&&\quad+\sum_{k\neq 0} \frac{i}{k\Omega}\partial_xf_0(X(t),X(t-\tau))\:f_k(X(t),X(t-\tau))\\
&&\quad-\sum_{k\neq 0} \frac{i}{k\Omega}\partial_xf_k(X(t),X(t-\tau))\:f_0(X(t),X(t-\tau))\\&&
\quad+\sum_{k\neq 0} \frac{i}{k\Omega}\partial_xf_k(X(t),X(t-\tau))\:f_{-k}(X(t),X(t-\tau))
\\&&
\quad-\sum_{k\neq 0} \frac{i}{k\Omega}\partial_yf_k(X(t),X(t-\tau))\:\dot\varphi(t-\tau),
\end{eqnarray*}
 and
\begin{eqnarray*}
F_{2,2}&=&f_0(X(t),X(t-\tau))\\
&&\quad+\sum_{k\neq 0} \frac{i}{k\Omega}\partial_xf_0(X(t),X(t-\tau))\:f_k(X(t),X(t-\tau))\\
&&\quad-\sum_{k\neq 0} \frac{i}{k\Omega}\partial_xf_k(X(t),X(t-\tau))\:f_0(X(t),X(t-\tau))\\&&
\quad+\sum_{k\neq 0} \frac{i}{k\Omega}\partial_xf_k(X(t),X(t-\tau))\:f_{-k}(X(t),X(t-\tau))
\\&&
\quad +\sum_{k\neq 0} \frac{i}{k\Omega}\partial_yf_0(X(t),X(t-\tau))\:f_{k}(X(t-\tau),X(t-2\tau))
\\&&
\quad -\sum_{k\neq 0} \frac{i}{k\Omega}\partial_yf_k(X(t),X(t-\tau))\:f_0(X(t-\tau),X(t-2\tau))
\\&&
\quad +\sum_{k\neq 0} \frac{i}{k\Omega}\partial_yf_k(X(t),X(t-\tau))\:f_{-k}(X(t-\tau),X(t-2\tau)),
\end{eqnarray*}
valid for \(t\geq \tau\). In these expressions \(\partial_x f_k(x,y)\) and \(\partial_y f_k(x,y)\)  denote the Jacobian matrices of
\(f_k(x,y)\) with respect to \(x\) and \(y\) respectively. By arguing as in the first order case, one  proves that at stroboscopic times
the difference between the solution \(x(t)\) of \eqref{eq:osc}--\eqref{eq:osc2} the solution \(X\) of \eqref{eq:aver12}--\eqref{eq:f22} is \(\mathcal{O}(1/\Omega^2)\).
\subsection{Third- and higher order averaged systems}
It is still feasible to obtain in closed form the third-order averaged system for delay problems starting
 from the corresponding  expression for the case without delay given in \cite{part2,erratum,icmat}. The
 expression of the averaged delay system one obtains in this way is lengthy and will not be reproduced here.
 In fact for order 3 or higher, the best approach is to find recursively the word basis functions
for the system  \eqref{eq:system}--\eqref{eq:system2} for the specific \(f\) of interest.

For an averaged system of order \(n\), the right-hand side of the system has \(n\) expressions corresponding to the intervals \(0\leq t<\tau\), \dots, \((n-2)\tau\leq t < (n-1)\tau\), and \(t\geq (n-1)\tau\). (Recall that the second-order averaged system displayed above has two expressions.)
\section{An example}
As an illustration, we consider the system
\begin{eqnarray*}
\frac{du}{dt}&=&\frac{\alpha}{1+v^{\beta}}-u(t-\tau)+A\sin(\omega t)+B\sin(\Omega t),\\
\frac{dv}{dt}&=&\frac{\alpha}{1+u^{\beta}}-v(t-\tau),
\end{eqnarray*}
where \(\alpha\) and \(\beta\) are parameters. There are two periodic forcing terms; the forcing \(A\sin(\omega t)\) is slow (i.e.\ the frequency \(\omega\) is of moderate size) and \(B\sin(\Omega t)\), \(\Omega \gg 1\), is fast. When there is no forcing (\(A=B=0\)), the system represents a delayed genetic toggle switch, a synthetic gene regulatory network \cite{gardner}. The forced system
was studied  in \cite{Daza} as an instance of the emergence of \emph{vibrational resonance} \cite{vr,guirao}, i.e. the enhancement, due to the presence of the fast forcing, of the response of the system to the slow forcing.

In order to have a system of the form \eqref{eq:osc} it is necessary to
rewrite \(A\sin(\omega t)\) as \(A\sin(\omega \hat t)\), where \(\hat t\) is a new dependent variable defined by the differential equation \((d/dt) \hat t  = 1\) and the initial condition \(\hat t = 0\) at  \(t=0\) (this is of course the technique used above in \eqref{eq:system2}--\eqref{eq:ic} to deal with the slow dependence on \(t\) of the right-hand side of \eqref{eq:system}). Thus only the fast forcing is averaged.
After computing recursively the required word basis functions and coefficients, the third-order averaged system is found to be given by
\begin{eqnarray*}
\frac{dU}{dt}&=&\frac{\alpha}{1+V^{\beta}}-U(t-\tau)+A\sin(\omega t),\\ \nonumber
\frac{dV}{dt}&=&\frac{\alpha}{1+U^{\beta}}-V(t-\tau)-\frac{B}{\Omega}\frac{\alpha\beta U^{\beta-1}}{(1+U^{\beta})^2}\\&&\qquad\qquad\qquad+\frac{B^2}{\Omega^2}\frac{3\alpha\beta U^{\beta-2}(U^{\beta}-\beta+\beta U^{\beta}+1)}{4(1+U^{\beta})^3},
\end{eqnarray*}
for $0\le t<\tau$, and
\begin{eqnarray*}\
\frac{dU}{dt}&=&\frac{\alpha}{1+V^{\beta}}-U(t-\tau)+A\sin(\omega t)-\frac{B}{\Omega},\\ \nonumber
\frac{dV}{dt}&=&\frac{\alpha}{1+U^{\beta}}-V(t-\tau)-\frac{B}{\Omega}\frac{\alpha\beta U^{\beta-1}}{(1+U^{\beta})^2}\\&&\qquad\qquad\qquad+\frac{B^2}{\Omega^2}\frac{3\alpha\beta U^{\beta-2}(U^{\beta}-\beta+\beta U^{\beta}+1)}{4(1+U^{\beta})^3},
\end{eqnarray*}
for $t\ge \tau$. (As pointed out above, in general, the third-order averaged system
has an analytic expression in \(\tau \leq t< 2\tau\) and a different analytic expression in
 \(t\geq 2\tau\); for the problem at hand those two expressions happen to coincide.) The second-order averaged system may be retrieved from the last displays by omitting the terms that have an \(\Omega^2\) factor in the denominator.
\begin{table}[t]
\footnotesize
\vspace{-0.5cm}
\begin{center}
\resizebox{\textwidth}{!}{
\begin{tabular}{rccccccccc}
\hline
 &$\Omega=16\pi$ & $\Omega=32\pi$ & $\Omega=64\pi$ & $\Omega=128\pi$ & $\Omega=256\pi$ & $\Omega=512\pi$
\\ \hline AS2&3.31(-4)&8.83(-5)& 2.28(-5) & 5.77(-6) & 1.46(-6) & 3.66(-7)\\ AS3&1.04(-4)&1.28(-5)& 1.59(-6) & 1.96(-7) & 2.07(-8) &2.30(-9) \\
\hline
\end{tabular}}
\end{center}
\label{newtable}
\vspace{-0.5cm}
\caption{\footnotesize
Maximum errors in $u$, with respect to the true solution, for the second-order and the third-order averaged solutions in the interval $0\leq t\leq 2$ }
\end{table}

We have integrated numerically  the true dynamics and the  second- and third-order averaged systems. The integrations were carried out with the Matlab function dde23 with relative and absolute tolerances $10^{-8}$ and $10^{-10}$ respectively.
For the choice $\alpha=2.5$, $\beta=2$,
$A=0.1$, $\omega=0.1$, $B=2.0$, and $\tau=0.5$, with constant history  \( u(t) = 0.5\), \(v(t) = 2.0\), \(-\tau \leq t\leq 0\) (an equilibrium of the unforced system),
Table~\ref{newtable} presents the maximum error in \(u\) at stroboscopic times in the short interval \(0\leq t\leq 2\); in agreement with the theory, the error behavior is \(\mathcal{O}(1/\Omega^2)\) for second-order averaging and \(\mathcal{O}(1/\Omega^3)\) for third-order averaging.

\begin{figure}[t]
\vspace{-4cm}\centering\includegraphics[scale=0.5]{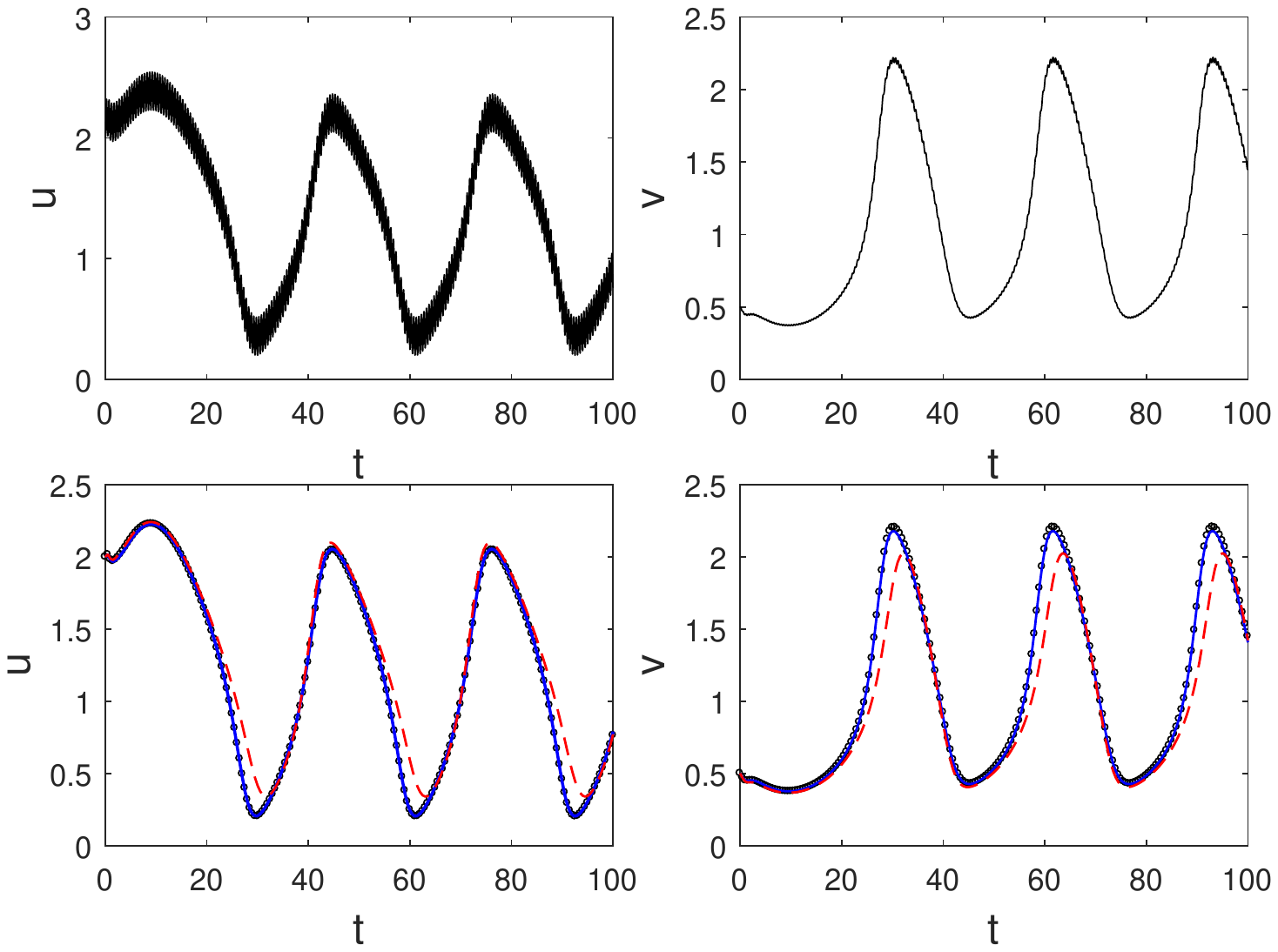}
\vspace{-4cm}
\caption{\footnotesize The left (respectively right) panels correspond to the \(u\) (respectively \(v\)) component of the solution. The top panels plot the true oscillatory solution as a function of \(t\). Due to the fast oscillations of large amplitude, the graph for \(u\)  appears as a solid band. The oscillations in \(v\) have a smaller amplitude; the differential equation for \(v\) does no include any periodic forcing and therefore the fast oscillations in \(v\) are only due to its coupling to \(u\). The bottom panels give the second order averaged solution (discontinuous line) and third order averaged solution (solid line). In the bottom panels, the true solution (circles) is represented only at stroboscopic times. Clearly the third order system  reproduces accurately the behaviour of the true solution (the circles are on the solid line). That is not the case if averaging is only carried out to second order.
}
\label{fig:C}
\end{figure}

In Fig.~\ref{fig:C}, the parameters are $\alpha=2.5$, $\beta=2$,
$A=0.2$, $\omega=0.2$, $B=2.0$, $\Omega=4\pi$ and $\tau=0.5$, with constant history
 \( u(t) = 2.0\), \(v(t) = 0.5\), \(-\tau \leq t\leq 0\) (again an equilibrium of the unforced system). The integration is perfomed in the interval \(0\leq t\leq 100\). The numerical integration of the oscillatory problem needed 261.0 seconds, a quantity that has to be compared with the 2.1 and 2.6 seconds required to integrate the second- and third-order averaged systems respectively. In studies like those performed in \cite{Daza}, where the oscillatory systems has to be integrated in long time intervals for many different choices of the values of the parameters, the advantage of averaging is then clear.

\section{Conclusion}
We have showed that, when the delay is a multiple of the forcing period, it is possible to extend to periodically forced, constant delay problems, the word-series approach to the systematic derivation of high-order averaged systems. We have presented an example where the new technique has been applied to a system of interest in connection with the phenomenon of vibrational resonance.

 To conclude let us point out that, for ordinary differential equations, it is possible to compute numerically a stroboscopically  averaged solution
\(\Xi\) without the explicit knowledge of the corresponding averaged system; the information on the averaged system required by the integrator is derived by numerically simulating the oscillatory system \eqref{eq:oscode} \cite{CCMS1,CCMS2}. Such techniques have been extended to delay differential equations \cite{beibei,mpbeibei}.

\bigskip

{\bf Acknowledgements.} J.M.S. was supported by project MTM2016-77660-P(AEI/ FEDER, UE) funded by MINECO (Spain).
 B. Z. is supported by the National Natural Science Foundation of China (Grant No. 11771438) and the Postdoctoral Fund Project of China (Grant No. 2018M641506) and in addition by the National Center for Mathematics and Interdisciplinary Sciences, CAS.

\end{document}